\documentclass[11pt]{amsart}

\theoremstyle{plain} \numberwithin{equation}{section}
\newtheorem{Theorem}{Theorem}
\newtheorem{Lemma}[Theorem]{Lemma}

\newtheorem{Proposition}[Theorem]{Proposition}

\theoremstyle{remark}

\title[Bari--Markus property]
{Bari--Markus property for Riesz projections of 1D periodic Dirac
operators}

\author{Plamen Djakov}

\author{Boris Mityagin}

\begin{document}

\address{Sabanci University, Orhanli,
34956 Tuzla, Istanbul, Turkey}

 \email{djakov@sabanciuniv.edu}

\address{Department of Mathematics,
The Ohio State University,
 231 West 18th Ave,
Columbus, OH 43210, USA}

\email{mityagin.1@osu.edu}

\begin{abstract}
The Dirac operators
$$ Ly = i \begin{pmatrix} 1
& 0 \\ 0 & -1
\end{pmatrix}
\frac{dy}{dx}  + v(x) y, \quad y = \begin{pmatrix} y_1\\y_2
\end{pmatrix},  \quad
x\in[0,\pi],$$ with  $L^2$-potentials
$$ v(x) = \begin{pmatrix} 0 & P(x) \\ Q(x) & 0 \end{pmatrix},
 \quad P,Q \in L^2 ([0,\pi]), $$ considered on $[0,\pi]$
with periodic, antiperiodic or Dirichlet boundary conditions $(bc)$,
have discrete spectra, and the Riesz projections
$$
S_N = \frac{1}{2\pi i} \int_{ |z|= N-\frac{1}{2} }  (z-L_{bc})^{-1}
dz, \quad P_n = \frac{1}{2\pi i} \int_{ |z-n|= \frac{1}{4} }
(z-L_{bc})^{-1} dz
$$
are well--defined for  $|n| \geq N$ if $N $ is sufficiently large.
It is proved that
$$\sum_{|n| > N} \|P_n - P_n^0\|^2 < \infty,
$$
where $P_n^0, \, n \in \mathbb{Z},$ are the Riesz projections of the
free operator.

Then, by the Bari--Markus criterion, the spectral Riesz
decompositions
 $$ f = S_N f +  \sum_{|n| >N} P_n f, \quad  \forall f \in L^2; $$
converge unconditionally  in $L^2.$

\end{abstract}


\maketitle

\section{Introduction}

The question for unconditional convergence of the spectral
decompositions is one of the central problems in Spectral Theory of
Differential Operators \cite{D58,DS3,Mar,Na1,Na2,Min}.

In the case of ordinary differential operators on a finite interval,
say $I = [0,\pi],$
\begin{equation}
\label{i1} \ell (y) = \frac{d^m y}{dx^m} + \sum_{k=0}^{m-2} q_k (x)
\frac{d^k y}{dx^k}, \quad q_k \in H^k (I),
\end{equation}
with strongly regular boundary conditions $(bc)$ the eigenfunction
decomposition
\begin{equation}
\label{i2} f(x) = \sum_k c_k (f) u_k (x), \quad    \ell (u_k) =
\lambda_k u_k, \;  u_k \in (bc),
\end{equation}
{\em converge  unconditionally} for every $f \in L^2 (I)$ (see
\cite{Mih,Kes,DS3}).

If $(bc)$ are regular but not strictly regular the system of root
functions (eigenfunctions and associated functions) in general is
not a basis in $L^2.$ But if the root functions are combined
properly in disjoint groups $B_n, \; \cup B_n = \mathbb{N},  $  then
the series
\begin{equation}
\label{i3} f(x) = \sum_n P_n f, \quad P_n  f = \sum_{k\in B_n} c_k
(f) u_k (x),
\end{equation}
{\em converges unconditionally in } $L^2 $  (see
\cite{Shk79,Shk83}).

Let us be more specific in the case of operators of second order
\begin{equation}
\label{i4} \ell (y) = y^{\prime \prime} +q(x) y, \quad 0 \leq x \leq
\pi.
\end{equation}
Then, Dirichlet $bc =Dir:\; y(0) = y(\pi) =0 $ is {\em  strictly
regular}; however, Periodic $ bc= Per^+: \;  y(0) = y(\pi), \;
y^\prime(0) = y^\prime(\pi) $ and Antiperiodic $ bc= Per^-: \;  y(0)
=- y(\pi), \; y^\prime(0) = - y^\prime(\pi) $ are {\em regular, but
not strictly regular.}

Analysis -- even if it becomes more difficult and technical -- could
be extended to singular potentials $q \in H^{-1}.$  A. Savchuk and
A. Shkalikov showed (\cite{SS03}, Theorems 2.7 and 2.8) that for
both Dirichlet bc or (properly understood) Periodic or Antiperiodic
$bc,$ the spectral decomposition (\ref{i3}) converges
unconditionally. An alternative proof of this result is given in
\cite{DM19}.

For Dirac operators (\ref{01}) the results on unconditional
convergence are sparse and not complete so far
\cite{Shk83,Mal,MalOr,TY01,TY02,HO06}.

The case of separate boundary conditions, at least for smooth
potential $v,$ has been studied in detail in \cite{Mal}. For
periodic (or antiperiodic) $bc \;$  B.Mityagin proved unconditional
convergence of the series (\ref{i3}) with $\dim P_n = 2, \; |n| \geq
N(v),$  for potentials $v \in H^b, \; b >1/2 $ -- see Theorem 8.8
\cite{Mit04} for a precise statement.

Our techniques from \cite{DM19} to analyze the resolvents
$(\lambda-L_{bc})^{-1}$ of Hill operators with the weakest (in
Sobolev scale) assumption $v \in H^{-1}$  on "smoothness" of the
potential are adjusted and extended in the present paper to Dirac
operators with potentials in $L^2.$ We prove (see Theorem \ref{thm1}
for a precise statement) that if $v \in L^2$ and $bc = Per^\pm, Dir
$ the sequence of deviations $\|P_n -P_n^0\|$  is in $\ell^2.$ Then,
the Bari--Markus criterion (see \cite{B,M}  or \cite{GK}, Ch.6,
Sect.5.3, Theorem 5.2)) shows that the spectral decomposition
\begin{equation}
\label{i5} f = S_N f +  \sum_{|n| >N} P_n f, \quad  \forall f \in
L^2,
\end{equation}
where, for $|n|  \geq  N(v),$
\begin{equation}
\label{i6} \dim P_n = \begin{cases}  2  &  bc = Per^{\pm}\\ 1  & bc
=Dir  \end{cases},
\end{equation}
{\em converge unconditionally.} This is Theorem \ref{thm2}, the main
result of the present paper.

Further analysis requires thorough discussion of the algebraic
structure of {\em regular}  and {\em strictly regular} $bc$ for
Dirac operators. Then we can claim a general statement which is an
analogue of (\ref{i5})--(\ref{i6}), or Theorem \ref{thm2}, with $bc=
Dir$ in case of strictly regular boundary conditions, and $bc
=Per^\pm $ in case of regular but not strictly regular boundary
conditions.  We will give all the details in another paper.

\section{Preliminary results}

Consider the Dirac operator on $ I= [0,\pi] $
\begin{equation}
\label{01} Ly = i
\begin{pmatrix}
1 & 0 \\ 0 & -1
\end{pmatrix}
\frac{dy}{dx}  + v(x)y,
\end{equation}
where
\begin{equation}
\label{02} v(x)=
\begin{pmatrix}
0  & P(x) \\ Q(x) & 0
\end{pmatrix},
\quad y =
\begin{pmatrix}
y_1\\ y_2
\end{pmatrix},
\end{equation}
and $v$ is an $L^2$--potential, i.e., $P,Q \in L^2 (I). $

We equip the space $ H^0 $ of $L^2 (I)$--vector functions $F=
\begin{pmatrix} f_1\\f_2 \end{pmatrix}$ with the scalar product $$
\langle F,G \rangle = \frac{1}{\pi} \int_0^\pi \left ( f_1(x)
\overline{g_1(x)} +  f_2(x) \overline{g_2(x)} \right ) dx. $$

Consider the following boundary conditions $(bc):$

(a) {\em periodic} $Per^+ : \quad y(0) = y(\pi),$ i.e., $ y_1 (0) =
y_1 (\pi) $ and $ y_2 (0) = y_2 (\pi); $

(b) {\em anti-periodic} $Per^- : \;  y(0) = - y(\pi),$ i.e., $ y_1
(0) = -y_1 (\pi) $ and $ y_2 (0) = -y_2 (\pi); $

(c) {\em Dirichlet} $Dir: \quad y_1 (0) = y_2 (0), \; y_1 (\pi) =
y_2 (\pi).$

The corresponding closed operator with a domain
\begin{equation}
\label{07}
 \Delta_{bc} = \left \{f \in (W_1^2 (I))^2  \; : \;
\; f =
\begin{pmatrix} f_1 \\ f_2 \end{pmatrix}
 \in (bc)  \right \}
\end{equation}
 will be denoted by $L_{bc},$ or respectively, by $L_{Per^\pm} $ and
$ L_{dir}.$ If $v=0,$ i.e., $P\equiv 0, Q \equiv 0,$ we write
$L^0_{bc}$ (or simply $L^0$), or $L^0_{Per^\pm}, L^0_{Dir}$
respectively. Of course, it is easy to describe the spectra and
eigenfunctions for $L^0_{bc}.$

(a)  $Sp (L^0_{Per^+}) = \{n \;\text{even} \} = 2 \mathbb{Z}; $ each
number $n \in 2\mathbb{Z} $ is a double eigenvalue, and the
corresponding eigenspace is
\begin{equation}
\label{011} E^0_n = Span \{e^1_n, e^2_n\},\quad n \in 2 \mathbb{Z},
\end{equation}
 where
\begin{equation}
\label{012}
 \displaystyle e^1_n (x) =
\begin{pmatrix}
e^{-inx}\\ 0
\end{pmatrix},
\quad e^2_n (x) =
\begin{pmatrix}
0\\ e^{inx}
\end{pmatrix};
\end{equation}

(b)  $Sp (L^0_{Per^-}) = \{n \;\text{odd} \} = 2 \mathbb{Z}+ 1; $
the corresponding eigenspaces $E_n^0 $ are given by (\ref{011}) and
(\ref{012}) but with $ n \in 2\mathbb{Z} +1; $

(c) $ Sp (L^0_{Dir}) = \{n \in \mathbb{Z} \};$ each eigenvalue $n $
is simple. The corresponding normalized eigenfunction is
\begin{equation}
\label{013} g_n  (x) = \frac{1}{\sqrt{2}} \left (e^1_n + e^2_n
\right ), \quad n \in \mathbb{Z},
\end{equation}
so the corresponding (one-dimensional) eigenspace is
\begin{equation}
\label{014} G_n^0 = Span \{ g_n \}.
\end{equation}
\vspace{3mm}

We study the spectral properties of the operators $L_{Per^\pm} $ and
$L_{Dir} $ by using their Fourier representations with respect to
the eigenvectors of the corresponding free operators given above in
(\ref{011})--(\ref{014}).

Let
\begin{equation}
\label{03} P(x) = \sum_{m \in 2\mathbb{Z}} p(m) e^{imx}, \qquad
 Q(x) = \sum_{m \in 2\mathbb{Z}} q(m) e^{imx},
\end{equation}
and
\begin{equation}
\label{04} P(x) = \sum_{m \in 1+2\mathbb{Z}} p_1(m) e^{imx}, \qquad
  Q(x) = \sum_{m \in 1+ 2\mathbb{Z}} q_1(m) e^{imx},
\end{equation}
 be, respectively, the Fourier expansions of the functions $P$ and $Q $
about the systems  $\{e^{imx}, \, m \in 2\mathbb{Z}\}$ and
$\{e^{imx}, \, m \in 1+2\mathbb{Z}\}.$

 Then
\begin{equation}
\label{05} \|v\|^2  =\sum_{m\in 2\mathbb{Z}} \left ( |p(m)|^2 +
|q(m)|^2 \right ) =\sum_{m\in 1+2\mathbb{Z}} \left ( |p_1(m)|^2 +
|q_1(m)|^2 \right ).
\end{equation}

In its Fourier representation, the operator $L^0 $  is diagonal, and
$V$ is defined by its action on vectors $e_n^1$ and $ e_n^2,$ with
$n \in 2\mathbb{Z}$ for $bc = Per^+$ and $n \in 1+ 2\mathbb{Z}$ for
$bc = Per^-.$ In view of (\ref{02}) and (\ref{03}), we have
\begin{equation}
\label{12} Ve_n^1 = \sum_{k  \in n + 2\mathbb{Z}} q(k+n) e_k^2,
\quad Ve_n^2 = \sum_{k  \in n + 2\mathbb{Z}} p(-k-n) e_k^1,
\end{equation}
so, the matrix representation of $V$ is
\begin{equation}
\label{13}  V \sim
\begin{pmatrix} 0 &  V^{12}\\ V^{21} & 0 \end{pmatrix}, \quad
 (V^{12})_{kn} = p(-k-n), \quad  (V^{21})_{kn} = q(k+n).
\end{equation}

In the case of Dirichlet boundary conditions the operator
 $L^0$ is diagonal as well. The matrix
representation of $V$ given by the following lemma.

\begin{Lemma}
\label{lem1} Let $(g_n)_{n \in \mathbb{Z}} $ be the orthogonal
normalized basis of eigenfunctions of $L^0 $ in the case of
Dirichlet boundary conditions. Then
\begin{equation}
\label{14} V_{kn}:= \langle Vg_n, g_k \rangle = W(k+n), \quad k,n
\in \mathbb{Z},
\end{equation}
with
\begin{equation}
\label{15} W(m) = \begin{cases} (p(-m) +q(m))/2     & m\; \text{even} \\
(p_1(-m) +q_1(m))/2  & m\; \text{odd}
\end{cases}.
\end{equation}
\end{Lemma}
The proof follows from a direct computation of $\langle Vg_n, g_k
\rangle.$   Let us mention, that the sequences  $p_1 (m)$ and $q_1
(m)$ in (\ref{15}) are Hilbert transforms of $p(n)$ and $q(n)$ (see
\cite{DM15}, Lemma 2 in Section 1.3) but we do not need this fact.
In the following only the relation (\ref{05}) is essential.

 In view of (\ref{011})--(\ref{014}) the operator $R^0_\lambda =
 (\lambda -L^0)^{-1} $
is well defined, respectively, for $\lambda \not \in 2\mathbb{Z}$ if
$bc = Per^+,$  $\lambda \not \in 1+ 2\mathbb{Z}$ if $bc = Per^-,$
and $\lambda \not \in \mathbb{Z}$ if $bc = Dir.$ The operator
$R^0_\lambda$ is diagonal, and we have
 \begin{equation}
\label{21} R^0_\lambda e^1_n = \frac{1}{\lambda -n} e_n^1, \quad
R^0_\lambda e^2_n = \frac{1}{\lambda -n} e_n^2 \quad \text{for} \;
bc= Per^\pm,
\end{equation}
and
 \begin{equation}
\label{22} R^0_\lambda g_n = \frac{1}{\lambda -n} g_n \quad
\text{for} \; bc = Dir.
\end{equation}

The standard perturbation type formulae for the resolvent $R_\lambda
= (\lambda - L^0 -V)^{-1}$ are
\begin{equation}
\label{23} R_\lambda = (1-R_\lambda^0 V)^{-1} R_\lambda^0 =
\sum_{k=0}^\infty (R_\lambda^0 V)^k R_\lambda^0,
\end{equation}
and
\begin{equation}
\label{24} R_\lambda = R_\lambda^0 (1-V R_\lambda^0 )^{-1}  =
\sum_{k=0}^\infty R_\lambda^0 (V R_\lambda^0)^k.
\end{equation}

 The simplest
conditions that guarantee convergence of the series (\ref{23}) or
(\ref{24}) in $\ell^2 $ are $$ \|R_\lambda^0 V\| <1, \quad
\mbox{respectively,} \quad \|V R_\lambda^0\| < 1. $$ In the case of
Dirac operators there are no such good estimates but there are good
estimates for the norms of $(R_\lambda^0 V)^2 $ and $(V
R_\lambda^0)^2 $ (see \cite{DM6} and \cite{DM15}, Section 1.2, for
more comments).

But now we are going to suggest another approach that is borrowed
from the study of Hill operators with periodic singular potentials
(see \cite{DM16,DM18,DM19}). Notice, that one can write (\ref{23})
or (\ref{24}) as
\begin{equation}
\label{25} R_\lambda = R^0_\lambda + R^0_\lambda V R^0_\lambda +
 \cdots = K^2_\lambda +
\sum_{m=1}^\infty K_\lambda(K_\lambda V K_\lambda)^m K_\lambda,
\end{equation}
provided
\begin{equation}
\label{26} (K_\lambda)^2 = R^0_\lambda.
\end{equation}

In view of (\ref{21}) and (\ref{22}), we define an operator $K=
K_\lambda $ with the property (\ref{26}) by
\begin{equation}
\label{31} K_\lambda e^1_n = \frac{1}{\sqrt{\lambda -n}} e_n^1,
\quad K_\lambda e^2_n = \frac{1}{\sqrt{\lambda -n}} e_n^2 \quad
\text{for} \; bc= Per^\pm,
\end{equation}
and
 \begin{equation}
\label{32} K_\lambda g_n = \frac{1}{\sqrt{\lambda -n}} g_n \quad
\text{for} \; bc = Dir,
\end{equation}
where $$\sqrt{z}= \sqrt{r} e^{i\varphi/2} \quad \mbox{if} \quad z=
re^{i\varphi}, \;\; -\pi \leq \varphi < \pi. $$ Then $R_\lambda $ is
well--defined if
\begin{equation}
\label{33} \|K_\lambda V K_\lambda \|_{\ell^2   \to \ell^2} <1.
\end{equation}

In view of (\ref{12}) and (\ref{31}), for periodic or anti--periodic
boundary conditions $bc = Per^\pm,$  we have
\begin{equation}
\label{34}  (K_\lambda V K_\lambda) e^1_n = \sum_k
\frac{q(k+n)}{(\lambda -k)^{1/2}(\lambda -n)^{1/2}} e_k^2, \quad
 (K_\lambda V K_\lambda) e^2_n = \sum_k
\frac{p(-k-n)}{(\lambda -k)^{1/2}(\lambda -n)^{1/2}} e_k^1,
\end{equation}
so, the Hilbert--Schmidt norm of the operator  $K_\lambda V
K_\lambda$ is given by
\begin{equation}
\label{35}  \|K_\lambda V K_\lambda\|^2_{HS}  = \sum_{k,m}
\frac{|q(k+m)|^2}{|\lambda -k||\lambda -m|} +\sum_{k,m}
\frac{|p(-k-m)|^2}{|\lambda -k||\lambda -m|},
\end{equation}
where $k,m \in 2\mathbb{Z}$ for $bc = Per^+$ and $k,m \in 1+
2\mathbb{Z}$ for $bc = Per^-.$

In an analogous way (\ref{14}), (\ref{15}) and (\ref{32}) imply, for
Dirichlet boundary conditions $bc =Dir,$
\begin{equation}
\label{37}  (K_\lambda V K_\lambda) g_n = \sum_k
\frac{W(k+n)}{(\lambda -k)^{1/2}(\lambda -n)^{1/2}} g_k, \quad k,n
\in \mathbb{Z},
\end{equation}
and therefore, we have
\begin{equation}
\label{38}  \|K_\lambda V K_\lambda\|^2_{HS}  = \sum_{k,m}
\frac{|W(k+m)|^2}{|\lambda -k||\lambda -m|},\quad k,m \in
\mathbb{Z}.
\end{equation}

For convenience, we set
\begin{equation}
\label{40} r(m) = \max (|p(m)|,|p(-m)|) + \max (|q(m)|,|q(-m)|),
\quad m \in 2\mathbb{Z},
\end{equation}
if $bc = Per^\pm,$ and
\begin{equation}
\label{41} r(m) =  |W(m)| \quad m \in \mathbb{Z},
\end{equation}
if $bc = Dir.$ Now we define operators $\bar{V} $ and
$\bar{K}_\lambda $ which dominate, respectively, $V$ and
$K_\lambda,$ as follows:
\begin{equation}
\label{42} \bar{V}e_n^1 = \sum_{k  \in n + 2\mathbb{Z}} r(k+n)
e_k^2, \quad \bar{V}e_n^2 = \sum_{k  \in n + 2\mathbb{Z}} r(k+n)
e_k^1 \quad \text{for} \; bc= Per^\pm,
\end{equation}
\begin{equation}
\label{43} \bar{V} g_n = \sum_{k  \in \mathbb{Z}} r(k+n) g_k \quad
\text{for} \; bc = Dir,
\end{equation}
and
\begin{equation}
\label{44} \bar{K}_\lambda e^1_n = \frac{1}{\sqrt{|\lambda -n|}}
e_n^1, \quad \bar{K}_\lambda e^2_n = \frac{1}{\sqrt{|\lambda -n|}}
e_n^2 \quad \text{for} \; bc= Per^\pm,
\end{equation}
 \begin{equation}
\label{45} \bar{K}_\lambda g_n = \frac{1}{\sqrt{|\lambda -n|}} g_n
\quad \text{for} \; bc = Dir.
\end{equation}

Since the matrix elements of the operator $K_\lambda V K_\lambda$ do
not exceed, by absolute value, the matrix elements of
$\bar{K}_\lambda \bar{V} \bar{K}_\lambda,$ we estimate from above
the Hilbert--Schmidt norm of the operator $K_\lambda V K_\lambda$ by
one and the same formula:
\begin{equation}
\label{46}
 \|K_\lambda V K_\lambda\|^2_{HS}  \leq
\|\bar{K}_\lambda \bar{V} \bar{K}_\lambda\|_{HS} =
 \sum_{i,k}
\frac{|r(i+k)|^2}{|\lambda -i||\lambda -k|},
\end{equation}
where  $i,k \in 2\mathbb{Z}$ if $bc = Per^+$ and  $i,k \in 1+
2\mathbb{Z}$ if $bc = Per^-,$ or $ i,k \in \mathbb{Z}$ if $bc =Dir.$
Next we estimate the Hilbert--Schmidt norm of the operator
$\bar{K}_\lambda \bar{V} \bar{K}_\lambda$ for  $ \lambda \in C_n =
\{\lambda:\; |\lambda -n|= 1/2 \}. $

For each $\ell^2$--sequence  $x=(x(j))_{j \in \mathbb{Z}} $ and $m
\in \mathbb{N}$ we set
\begin{equation}
\label{50} \mathcal{E}_m (x) = \left (   \sum_{|j|\geq m} |x(j)|^2
\right )^{1/2}.
\end{equation}

\begin{Lemma}
\label{lem2} In the above notations, if $n\neq 0,$ then
\begin{equation}
\label{51} \|\bar{K}_\lambda \bar{V} \bar{K}_\lambda\|^2_{HS} =
 \sum_{i,k}
\frac{|r(i+k)|^2}{|\lambda -i||\lambda -k|} \leq  C \left (
\frac{\|r\|^2}{\sqrt{|n|}} + (\mathcal{E}_{|n|} (r))^2 \right ),
\quad \lambda \in C_n,
\end{equation}
where $C$ is an absolute constant.
\end{Lemma}

{\em Remark:} For convenience, here and thereafter we denote by  $C$
any absolute constant.

\begin{proof}
Since
\begin{equation}
\label{52} 2|\lambda - i| \geq |n-i|  \quad \text{if} \; i \neq n,
\; \lambda \in C_n= \{\lambda:\; |\lambda -n|= 1/2 \},
\end{equation}
the sum in (\ref{51}) does not exceed
$$
4 |r(2n)|^2 + 4\sum_{k\neq n} \frac{|r(n+k)|^2}{|n-k|} +
4\sum_{i\neq n} \frac{|r(n+i)|^2}{|n-i|} + 4\sum_{i,k\neq n}
\frac{|r(i+k)|^2}{|n-i||n-k|}.
$$
In view of (\ref{t1}) and (\ref{t2}) in Lemma \ref{lemt1}, each of
the above sums does not exceed the right-hand side of (\ref{51}),
which completes the proof.
\end{proof}

{\em Corollary: There is $N \in \mathbb{N}$ such that}
\begin{equation}
\label{53} \|K_\lambda V K_\lambda\| \leq 1/2  \quad \text{for}\;\;
\lambda \in C_n, \; |n| >N.
\end{equation}

\section{Core results}

 By our Theorem 18 in \cite{DM15} (about spectra localization),
 for sufficiently large $|n|,$ say $|n|>N,$ the operator
$L_{Per^\pm}$ has exactly two (counted with their algebraic
multiplicity) periodic (for even $n$) or antiperiodic (for odd $n$)
eigenvalues  inside the disc with a center $n$ of radius $1/2.$ The
operator $L_{Dir}$ has, for all sufficiently large $|n|,$ one
eigenvalue in every such disc.

Let $P_n$ and $P_n^0 $ be the Riesz projections corresponding to $L$
and $L^0,$ i.e.,
$$ P_n = \frac{1}{2\pi i} \int_{C_n} (\lambda - L)^{-1} d\lambda,
\quad P_n^0 = \frac{1}{2\pi i} \int_{C_n} (\lambda - L)^{-1}
d\lambda,
$$ where $C_n = \{\lambda: \; |\lambda - n| =1/2 \}.   $

\begin{Theorem}
\label{thm1} Suppose $L$ and $L^0$ are, respectively, the Dirac
operator (\ref{01}) with an $L^2$ potential $v$ and the free Dirac
operator, subject to periodic, antiperiodic or Dirichlet boundary
conditions $bc =Per^\pm $ or $Dir.$ Then, there is $N \in
\mathbb{N}$ such that for $|n| >N $ the Riesz projections $P_n$ and
$P_n^0 $ corresponding to $L$ and $L^0$ are well defined and we have
\begin{equation}
\label{110} \sum_{|n|>N} \|P_n - P_n^0\|^2 < \infty.
\end{equation}
\end{Theorem}

\begin{proof}
Now we present the proof of the theorem up to a few technical
inequalities. They will be proved  later in Section 4, Lemmas
\ref{lemt1} and \ref{lemt2}. \vspace{3mm}

1. Let us notice that the operator--valued function $K_\lambda$ is
analytic in $ \mathbb{C} \setminus \mathbb{R}_+$. But (\ref{25}),
(\ref{111}) below and all formulas of this section -- which are
essentially variations of (\ref{25}) -- always have even powers of
$K_\lambda,$ and $K_\lambda^2 = R_\lambda^0 $ is analytic on
$\mathbb{C}\setminus Sp(L^0).$ Certainly, this justifies the use of
Cauchy formula or Cauchy theorem when warranted.

In view of (\ref{53}), the corollary after the proof of Lemma 2, if
$|n|$ is sufficiently large then the series in (\ref{25}) converges.
Therefore,
\begin{equation}
\label{111} P_n - P_n^0 = \frac{1}{2\pi i} \int_{C_n}
\sum_{s=0}^\infty K_\lambda   (K_\lambda V  K_\lambda )^{s+1}
K_\lambda d\lambda.
\end{equation}

{\em Remark.} We are going to prove (\ref{110}) by estimating the
Hilbert--Schmidt norms $\|P_n - P_n^0\|_{HS}$  which dominate $\|P_n
- P_n^0\|.$ Of course, these norms are equivalent as long as the
dimensions $dim \,(P_n - P_n^0)$ are uniformly bounded because for
any finite dimensional operator $T $ we have
$$
\|T\| \leq \|T\|_{HS} \leq (dim \,T)^{1/2} \|T\|
$$
but in the context of this paper  for all projections $dim \,P_n,
\;dim \,P_n^0 \leq 2.$
 \vspace{3mm}

2.  If $bc = Dir,$ then, by (\ref{013}),
$$
\|P_n - P_n^0\|_{HS}^2 = \sum_{m,k \in \mathbb{Z}} |\langle (P_n -
P_n^0)g_m,g_k \rangle|^2.
$$
By (\ref{111}), we get
$$
\langle (P_n - P_n^0)g_m,g_k \rangle = \sum_{s=0}^\infty I_n(s,k,m),
$$
where
$$
I_n(s,k,m)= \frac{1}{2\pi i} \int_{C_n}\langle K_\lambda(K_\lambda V
K_\lambda )^{s+1}K_\lambda g_m, g_k \rangle d\lambda.
$$
Therefore,
$$
\sum_{|n|>N} \|P_n - P_n^0\|_{HS}^2 \leq  \sum_{s,t=0}^\infty
\sum_{|n|>N} \sum_{m,k \in \mathbb{Z}} |I_n(s,k,m)|\cdot
|I_n(t,k,m)|.
$$
Now, the Cauchy inequality implies
\begin{equation}
\label{112} \sum_{|n|>N} \|P_n - P_n^0\|_{HS}^2 \leq
 \sum_{s,t=0}^\infty (A(s))^{1/2}(A(t))^{1/2},
\end{equation}
where
\begin{equation}
\label{113} A(s)=\sum_{|n|>N}\sum_{m,k \in \mathbb{Z}}
|I_n(s,k,m)|^2.
\end{equation}

Notice that $A(s)$ depends on $N$ but this dependence is suppressed
in the notation.

 From the matrix
representation of the operators $K_\lambda $ and $V$ we get
\begin{equation}
\label{128} \langle K_\lambda (K_\lambda V K_\lambda )^{s+1}
K_\lambda e_m, e_k \rangle
 =  \sum_{j_1,\ldots,j_s} \frac{W(k + j_1) W(j_1+j_2) \cdots W(j_s
+m)}{(\lambda -k)(\lambda -j_1)\cdots (\lambda -j_s)(\lambda -m)},
\end{equation}
and therefore,
\begin{equation}
\label{129} I_n (s,k,m)=    \frac{1}{2\pi i}\int_{C_n} \sum_{j_1,
\ldots j_s} \frac{W(k + j_1) W(j_1+j_2) \cdots W(j_s +m)}{(\lambda
-k)(\lambda -j_1)\cdots (\lambda -j_s)(\lambda -m)} d\lambda.
\end{equation}
In view of (\ref{41}), we have
\begin{equation}
\label{130} \left |\frac{W(k + j_1) W(j_1+j_2) \cdots W(j_s
+m)}{(\lambda -k)(\lambda -j_1)\cdots (\lambda -j_s)(\lambda -m)}
\right | \leq   B(\lambda, k,j_1,\ldots,j_s,m),
\end{equation}
where
\begin{equation}
\label{133} B(\lambda, k,j_1,\ldots,j_s,m)= \frac{r(k+j_1)r(j_1
+j_2) \cdots r(j_{s-1}+j_s) r(j_s +m)}{|\lambda -k||\lambda -j_1|
\cdots |\lambda -j_s| |\lambda -m |},   \quad s>0,
\end{equation}
and
\begin{equation}
\label{134} B(\lambda,k,m)= \frac{r(m+k)}{|\lambda -k||\lambda -m |}
\end{equation}
in the case when $s=0 $ and there are no $j$-indices.
 Moreover, by (\ref{41}),(\ref{43})  and (\ref{45}), we have
\begin{equation}
\label{135} \sum_{j_1,\ldots,j_s} B(\lambda,k, j_1,\ldots,j_s,m) =
\langle \bar{K}_\lambda (\bar{K}_\lambda \bar{V}
\bar{K}_\lambda)^{s+1}\bar{K}_z e_m, e_k \rangle.
\end{equation}

\begin{Lemma}
\label{lemd} In the above notations, we have
\begin{equation}
\label{136} A(s) \leq  B_1 (s)+B_2 (s)+B_3 (s)+B_4 (s),
\end{equation}
where
\begin{equation}
\label{137} B_1 (s)= \sum_{|n|>N} \sup_{\lambda \in C_n} \left (
\sum_{j_1, \ldots, j_s}
 B(\lambda, n,j_1,\ldots,j_s,n) \right )^2;
\end{equation}
\begin{equation}
\label{138} B_2 (s)= \sum_{|n|>N} \sum_{k\neq n}\sup_{\lambda \in
C_n}\left ( \sum_{j_1, \ldots, j_s}  B(\lambda, k,j_1,\ldots,j_s,n)
\right )^2;
\end{equation}
\begin{equation}
\label{139} B_3 (s)= \sum_{|n|>N} \sum_{m\neq n} \sup_{\lambda \in
C_n} \left ( \sum_{j_1, \ldots, j_s}  B(\lambda, n,j_1,\ldots,j_s,m)
\right )^2;
\end{equation}
\begin{equation}
\label{140}  B_4 (s)= \sum_{|n|>N} \sum_{m,k\neq n}\sup_{\lambda \in
C_n}\left (\sum_{j_1, \ldots, j_s}^*  B(\lambda, k,j_1,\ldots,j_s,m)
\right )^2,   \quad s\geq 1,
\end{equation}
where the symbol $*$ over the sum in the parentheses means that at
least one of the indices $j_1, \ldots, j_s $ is equal to $n.$
\end{Lemma}

\begin{proof}
Indeed, in view of (\ref{113}), we have
$$  A(s)  \leq  A_1 (s)+A_2 (s)+A_3 (s)+A_4 (s),$$
where
$$ A_1 (s) = \sum_{|n|>N} |I_n (s,n,n)|^2, \quad
A_2 (s) = \sum_{|n|>N} \sum_{k\neq n} |I_n (s,k,n)|^2 $$
$$A_3 (s) = \sum_{|n|>N} \sum_{m\neq n} |I_n (s,n,m)|^2, \quad
    A_4 (s) = \sum_{|n|>N} \sum_{m,k\neq n} |I_n (s,k,m)|^2.     $$
By (\ref{129})--(\ref{134}) we get immediately that
$$ A_\nu (s) \leq B_\nu (s), \quad  \nu= 1,2,3.$$
On the other hand,  by the Cauchy formula,
$$
\int_{C_n}  \frac{W(k + j_1) W(j_1+j_2) \cdots W(j_s +m)}{(\lambda
-k)(\lambda -j_1)\cdots (\lambda -j_s)(\lambda -m)} d\lambda = 0
\quad \text{if} \quad k, j_1, \ldots, j_s, m \neq n.
$$
Therefore, removing  from the sum in (\ref{129}) the terms with zero
integrals, and estimating from above the remaining sum, we get
$$ |I_n (s,k,m)| \leq \sup_{\lambda \in
C_n}\left (\sum_{j_1, \ldots, j_s}^*  B(\lambda, k,j_1,\ldots,j_s,m)
\right ), \qquad m,k \neq n. $$ From here it follows that $A_4 (s)
\leq B_4 (s),$ which completes the proof.
\end{proof}

3. In view of (\ref{112}) and (\ref{136}), Theorem \ref{thm1} will
be proved if we get ''good estimates'' of the sums $B_{\nu} (s), \;
\nu =1,\ldots,4, $ that are defined by (\ref{137})--(\ref{140}).

If $bc = Per^\pm,$ then using the orthonormal system of eigenvectors
of the free operator $L^0$ given by  (\ref{012}), we get
\begin{equation}
\label{0111} \|P_n - P_n^0\|_{HS}^2 =\sum_{\alpha, \beta=1}^2
\sum_{m,k} |\langle (P_n - P_n^0)e^\alpha_m,e^\beta_k \rangle|^2,
\end{equation}
where $m,k \in 2\mathbb{Z}$ if $n$ is even or $m,k \in 1+
2\mathbb{Z}$ if $n$ is odd. By (\ref{111}), we have
\begin{equation}
\label{0112} \langle (P_n - P_n^0)e^\alpha_m,e^\beta_k \rangle =
\sum_{s=0}^\infty I^{\alpha \beta}(n,s,k,m),
\end{equation}
where
\begin{equation}
\label{0113} I^{\alpha \beta}(n,s,k,m)= \frac{1}{2\pi i}
\int_{C_n}\langle K_\lambda(K_\lambda V K_\lambda )^{s+1}K_\lambda
e^\alpha_m,e^\beta_k \rangle d\lambda.
\end{equation}
Therefore,
$$
\sum_{|n|>N} \|P_n - P_n^0\|_{HS}^2 \leq  \sum_{\alpha, \beta=1}^2
\sum_{t,s=0}^\infty \sum_{|n|>N} \sum_{m,k} |I^{\alpha
\beta}(n,s,k,m)|\cdot |I^{\alpha \beta}(n,t,k,m)|.
$$
Now, the Cauchy inequality implies
\begin{equation}
\label{0114} \sum_{|n|>N} \|P_n - P_n^0\|_{HS}^2 \leq \sum_{\alpha,
\beta=1}^2
 \sum_{t,s=0}^\infty (A^{\alpha \beta}(s))^{1/2}(A^{\alpha \beta}(t))^{1/2},
\end{equation}
where
\begin{equation}
\label{141} A^{\alpha \beta}(s)=\sum_{|n|>N}\sum_{m,k} |I^{\alpha
\beta}(n,s,k,m)|^2.
\end{equation}

\begin{Lemma}
\label{lemp} In the above notations, with $r$ given by (\ref{40}),
$B(\lambda, k,j_1,\ldots,j_s,m)$ defined in (\ref{133}),(\ref{134}),
and $B_j (s), \, j=1, \ldots, 4,$ defined by
(\ref{137})--(\ref{140}), we have
\begin{equation}
\label{142} A^{\alpha \beta}(s) \leq  B_1 (s)+B_2 (s)+B_3 (s)+B_4
(s), \quad \alpha, \beta = 1,2.
\end{equation}
\end{Lemma}

\begin{proof}
The matrix representations of the operators $V$ and $K_\lambda  $
given in (\ref{13}) and (\ref{31}) imply that if $s $ is even, then
$\langle \langle K_\lambda(K_\lambda V K_\lambda )^{s+1}K_\lambda
e^\alpha_m,e^\alpha_k \rangle=0 $ for $ \alpha =1,2, $ and if $s$ is
odd then
\begin{equation}
\label{143} \langle K_\lambda(K_\lambda V K_\lambda )^{s+1}K_\lambda
e^1_m,e^1_k \rangle =\sum_{j_1,\ldots,j_s}\frac{p(-k-j_1 ) q(j_1 +
j_2 )\cdots p(-j_{s-1}- j_s ) q(j_s + m )}{(\lambda-k)(\lambda-j_1)
\cdots (\lambda-j_s)(\lambda-m) },
\end{equation}
\begin{equation}
\label{144} \langle K_\lambda(K_\lambda V K_\lambda )^{s+1}K_\lambda
e^2_m,e^2_k \rangle =\sum_{j_1,\ldots,j_s}\frac{q(k+j_1 ) p(-j_1 -
j_2 )\cdots q(j_{s-1}+ j_s ) p(-j_s - m )}{(\lambda-k)(\lambda-j_1)
\cdots (\lambda-j_s)(\lambda-m) }.
\end{equation}
In analogous way it follows that if $s$ is odd then $\langle \langle
K_\lambda(K_\lambda V K_\lambda )^{s+1}K_\lambda e^1_m,e^2_k
\rangle=0 $ and  $\langle \langle K_\lambda(K_\lambda V K_\lambda
)^{s+1}K_\lambda e^2_m,e^1_k \rangle=0, $ and if $s$ is even then
\begin{equation}
\label{145} \langle K_\lambda(K_\lambda V K_\lambda )^{s+1}K_\lambda
e^1_m,e^2_k \rangle =\sum_{j_1,\ldots,j_s}\frac{q(k+j_1 ) p(-j_1
-j_2 )\cdots p(-j_{s-1}- j_s ) q(j_s + m )}{(\lambda-k)(\lambda-j_1)
\cdots (\lambda-j_s)(\lambda-m) },
\end{equation}
\begin{equation}
\label{146} \langle K_\lambda(K_\lambda V K_\lambda )^{s+1}K_\lambda
e^2_m,e^1_k \rangle =\sum_{j_1,\ldots,j_s}\frac{p(-k-j_1 ) q(j_1
+j_2 )\cdots q(j_{s-1}+ j_s ) p(-j_s - m )}{(\lambda-k)(\lambda-j_1)
\cdots (\lambda-j_s)(\lambda-m) }.
\end{equation}
From (\ref{40}), (\ref{137})-(\ref{140}) and the above formulas it
follows that
$$
|\langle K_\lambda(K_\lambda V K_\lambda )^{s+1}K_\lambda
e^\alpha_m,e^\beta _k \rangle| \leq \sum_{j_1,\ldots,j_s} B(\lambda,
k,j_1,\ldots,j_s,m),$$ which implies immediately
\begin{equation}
\label{0146}  |I^{\alpha \beta}_n (s,k,m)| \leq \sup_{\lambda \in
C_n}\left (\sum_{j_1, \ldots, j_s}  B(\lambda, k,j_1,\ldots,j_s,m)
\right ).
\end{equation}

By (\ref{141}), $$  A^{\alpha \beta}(s)  \leq  A^{\alpha \beta}_1
(s)+A^{\alpha \beta}_2 (s)+A^{\alpha \beta}_3 (s)+A^{\alpha \beta}_4
(s),$$ where
$$ A^{\alpha \beta}_1 (s) = \sum_{|n|>N} |I^{\alpha \beta}_n (s,n,n)|^2, \quad
A^{\alpha \beta}_2 (s) = \sum_{|n|>N} \sum_{k\neq n} |I^{\alpha
\beta}_n (s,k,n)|^2
$$
$$A^{\alpha \beta}_3 (s) = \sum_{|n|>N} \sum_{m\neq n}
|I^{\alpha \beta}_n (s,n,m)|^2, \quad
    A^{\alpha \beta}_4 (s) = \sum_{|n|>N} \sum_{m,k\neq n}
    |I^{\alpha \beta}_n (s,k,m)|^2.     $$
Therefore, in view of (\ref{0146}) and (\ref{137})--(\ref{139}), we
get
$$ A^{\alpha \beta}_\nu (s) \leq B_\nu (s), \quad \nu =1,2,3. $$
Finally, as in the proof of Lemma~\ref{lemd}, we take into account
that in the sums (\ref{143})--(\ref{146}) the terms with indices
$j_1, \ldots, j_s, m,k \neq n  $ have zero integrals over the
contour $C_n. $ Therefore, $$ |I^{\alpha \beta}_n (s,k,m)| \leq
\sup_{\lambda \in C_n}\left (\sum^*_{j_1, \ldots, j_s}  B(\lambda,
k,j_1,\ldots,j_s,m) \right ), \quad m,k \neq n. $$  In view of
(\ref{140}), this yields $ A^{\alpha \beta}_4 (s) \leq B_4 (s),$
which completes the proof.
\end{proof}

Such estimates are given in the next proposition. For convenience,
we set for any $\ell^2 $--sequence $r= (r(j))$
\begin{equation}
\label{147} \rho_N = \left ( \frac{\|r\|^2}{\sqrt{N}} +
(\mathcal{E}_{N} (r))^2 \right )^{1/2}.
\end{equation}

\begin{Proposition}
\label{prop1} In the above notations,
\begin{equation}
\label{148}  B_\nu (s) \leq      C\|r\|^2 \rho_N^{2s},\quad \nu
=1,2,3, \qquad  B_4 (s) \leq  C s \|r\|^4 \rho_N^{2(s-1)}, \; s\geq
1,
\end{equation}
where $C$ is an absolute constant.
\end{Proposition}

\begin{proof}

{\em Estimates for $B_1 (s).$} By (\ref{134}) and (\ref{137}), we
have
$$ B_1 (0)= \sum_{|n|>N} \sup_{\lambda \in C_n}
\frac{|r(2n)|^2}{|\lambda -n|^2} = 4 (\mathcal{E}_{N} (r))^2 \leq
4\|r\|^2,
$$ so (\ref{148}) holds for $ B_1 (s)$ if $ s=0.$

If $s=1,$ then by (\ref{133}), the sum $B_1 (1)$ from (\ref{137})
has the form
$$B_1 (1) = \sum_{|n|>N} \sup_{\lambda \in C_n} \left | \sum_j
\frac{r(n+j)r(j+n)}{|\lambda -n||\lambda - j||\lambda -n|}  \right
|^2.
$$

By (\ref{52}), and  since $|\lambda -n| = 1/2 $ for $\lambda \in
C_n,$ we get
$$
B_1 (1) \leq \sum_{|n|>N}  \left (8\sum_{j\neq  n}
\frac{|r(j+n)|^2}{ |j-n|} +  8|r(2n)|^2\right )^2 $$
$$ \leq 128 \sum_{|n|>N}\left (
\sum_{j\neq  n} \frac{|r(j+n)|^2}{ |j-n|} \right )^2  +
128\sum_{|n|>N}|r(2n)|^4. $$  By the Cauchy inequality and
(\ref{t11}) in Lemma \ref{lemt2}, we have
$$
\sum_{|n|>N}\left ( \sum_{j\neq  n} \frac{|r(j+n)|^2}{ |j-n|} \right
)^2 \leq \sum_{|n|>N} \sum_{j\neq n} \frac{|r(j+n)|^2}{ |j-n|^2}
\|r\|^2 \leq C \|r\|^2\rho_N^2.
 $$
On the other hand, $\sum_{|n|>N}|r(2n)|^4 \leq \|r\|^2
(\mathcal{E}_N (r))^2 \leq \|r\|^2\rho_N^2, $ so (\ref{148}) holds
for $ B_1 (s)$ if $ s=1.$

Next, we consider the case  $s>1.$ In view of (\ref{133}), since
$|\lambda -n|= 1/2$ for $\lambda \in C_n, $
 the sum $B_1 (s)$ from (\ref{137}) can be written
 as
$$
B_1 (s) = \sum_{|n|>N} 4 \sup_{\lambda \in C_n} \left |\sum_{j_1,
\ldots,j_s} \frac{r(n+j_1)r(j_1+j_2)\cdots r(j_s
+n)}{|\lambda-j_1||\lambda-j_2|\cdots |\lambda-j_s|}\right |^2.
$$
Therefore, we have (with $j=j_1, k=j_s$)
$$
B_1 (s) = 4\sum_{|n|>N}  \sup_{\lambda \in C_n} \left |\sum_{j,k}
\frac{r(n+j)}{|\lambda - j|^{1/2}}\cdot H_{jk} (\lambda) \cdot
\frac{r(k +n)}{|\lambda - k|^{1/2}} \right |^2,
$$
where $ (H_{jk}(\lambda))$ is the matrix representation of the
operator $ H(\lambda) = (\bar{K}_\lambda \bar{V}
\bar{K}_\lambda)^{s-1}. $ By (\ref{51}) in Lemma \ref{lem2},
$$
\|H(\lambda)\|_{HS} =\left ( \sum_{j,k} |H_{jk}(\lambda)|^2 \right
)^{1/2} \leq \|\bar{K}_\lambda \bar{V} \bar{K}_\lambda \|_{HS}^{s-1}
\leq \rho_N^{s-1} \quad  \text{for} \; \lambda \in C_n, \; |n|>N. $$
Therefore, the Cauchy inequality implies
$$
B_1 (s) \leq  4\|H(\lambda)\|^2_{HS} \cdot \sigma \leq
4\rho_N^{2(s-1)}\cdot \sigma,
$$
where
$$\sigma=
\sum_{|n|>N} \sup_{\lambda \in C_n} \sum_{j,k}
\frac{|r(n+j)|^2}{|\lambda - j|}\cdot \frac{|r(k+n)|^2}{|\lambda -
k|}.
$$
By (\ref{52}) and since $|\lambda -n| = 1/2 $ for $\lambda \in C_n,$
we have
$$ \sigma \leq
  4\sum_{|n|>N}\sum_{j,k \neq n}
\frac{|r(n+j)|^2|r(n+k)|^2}{|n - j||n-k|}+
4\sum_{|n|>N}|r(2n)|^2\sum_{k \neq n} \frac{|r(n+k)|^2}{|n-k|}
$$
$$
 +
4\sum_{|n|>N}|r(2n)|^2\sum_{j \neq n} \frac{|r(n+j)|^2}{|n - j|} +
4\sum_{|n|>N}|r(2n)|^4.
$$
In view of (\ref{t12}) in Lemma \ref{lemt2}, the triple sum does not
exceed $C\|r\|^2\rho_N^2.$ By (\ref{t1}) in Lemma \ref{lemt1}, each
of the double sums can be estimated from above by
$$C\sum_{|n|>N}|r(2n)|^2 \rho_N^2 \leq C\|r\|^2\rho_N^2,  $$
and the same estimate holds for the single sum. Therefore,
$$
B_1 (s)  \leq C\rho_N^{2(s-1)}\cdot \|r\|^2 \rho_N^2,
$$
which completes the proof of (\ref{148}) for $B_1 (s).$\vspace{3mm}

{\em Estimates for $B_2 (s).$} By (\ref{134}) and (\ref{137}), we
have
$$ B_2 (0)= \sum_{|n|>N}\sum_{k \neq  n} \sup_{\lambda \in C_n}
\frac{|r(k+n)|^2}{|\lambda -k|^2|\lambda -n|^2}.
$$
Taking into account that  $|\lambda -n| = 1/2 $ for $\lambda \in
C_n,$ we get, in view of (\ref{52}) and (\ref{t11}) in Lemma
\ref{lemt2},
$$
B_2 (0) \leq 4\sum_{|n|>N}\sum_{k \neq  n}
\frac{|r(k+n)|^2}{|n-k|^2} \leq C \|r\|^2.
$$
So, (\ref{148}) holds for $B_2 (s)$ if $ s=0.$

If $s=1,$  then, by (\ref{133}), the sum $B_2 (s) $ in (\ref{148})
has the form
$$B_2 (1)=
\sum_{|n|>N} \sum_{ k \neq n} \sup_{\lambda \in C_n} \left |\sum_j
\frac{r(k+j)r(j+n)}{|\lambda -k||\lambda -j||\lambda -n|} \right
|^2.
$$
Since  $|\lambda -n| = 1/2 $ for $\lambda \in C_n,$ we get, in view
of (\ref{52}),
$$
B_2 (1) \leq \sum_{|n|>N} \sum_{ k \neq n} \left |\sum_{j\neq n} 8
\frac{r(k+j)r(j+n)}{|n -k||n -j|}   + 8 r(2n) \frac{r(k+n)}{|n-k|}
\right |^2.
$$
Therefore,
$$ B_2 (1) \leq 128\sigma_1 + 128\sigma_2,
$$
where (by the Cauchy inequality and (\ref{t11}) in Lemma
\ref{lemt2})
$$
\sigma_1 = \sum_{|n|>N,  k \neq n} \left (\sum_{j\neq n}
\frac{r(k+j)r(j+n)}{|n -k||n -j|} \right )^2 \leq \sum_{|n|>N, k\neq
n} \frac{1}{|n-k|^2}\left (\sum_{j\neq n} \frac{|r(n+j)|^2}{|n
-j|^2} \right ) \cdot \|r\|^2
$$
$$
=\sum_{|n|>N,  j \neq n} \frac{|r(n+j)|^2}{|n -j|^2} \sum_{k\neq n}
\frac{\|r\|^2}{|n-k|^2} \leq C \rho_N^2 \|r\|^2,
$$
and
$$
 \sigma_2 = \sum_{|n|>N,  k \neq n} |r(2n)|^2
\frac{|r(n+k)|^2}{|n-k|^2}\leq C \rho_N^2 \|r\|^2.
$$
Thus, (\ref{148}) holds for $B_2 (s) $ if $ s=1.$

If $s>1,$  then by (\ref{133}) and $|\lambda -n| =1/2$ for $\lambda
\in C_n,$ we have
$$B_2 (s) =
\sum_{|n|>N,  k \neq n} 2\sup_{\lambda \in C_n} \left
|\sum_{j_1,\ldots,j_s} \frac{r(k+j_1)r(j_1+j_2) \cdots r(j_s
+n)}{|\lambda - k||\lambda -j_1||\lambda -j_2| \cdots |\lambda -j_s|
} \right |^2.
$$
In view of (\ref{43}) and (\ref{44}), we get (with $j=j_1, i=j_s$)
$$B_2 (s) =
2\sum_{|n|>N,  k \neq n} \sup_{\lambda \in C_n} \left |\sum_{j,i}
\frac{r(k+j)}{|\lambda -k ||\lambda -j|^{1/2} } \cdot H_{ji}
(\lambda) \cdot \frac{r(i +n)}{|\lambda -i|^{1/2}} \right |^2,
$$
where $H_{ji} (\lambda)$ is the matrix representation of the
operator $H(\lambda) = (\bar{K}_\lambda \bar{V}
\bar{K}_\lambda)^{s-1}.$ Therefore, by the Cauchy inequality and
(\ref{51}) in Lemma \ref{lem2},
\begin{equation}
\label{165} B_2 (s) \leq 2 \sup_{\lambda \in
C_n}\|H(\lambda)\|^2_{HS} \cdot \tilde{\sigma}\leq 2 \rho_N^{2(s-1)}
\cdot \tilde{\sigma},
\end{equation}
where
$$
\tilde{\sigma} = \sum_{|n|>N,  k \neq n} \sup_{\lambda \in C_n}
\sum_{i,j} \frac{|r(k+j)|^2 |r(i+n)|^2}{|\lambda - k|^2|\lambda -
j||\lambda - i|}.
$$
From $|\lambda -n| = 1/2 $ for $\lambda \in C_n$ and (\ref{52}) it
follows that
$$
\tilde{\sigma} \leq 8(\tilde{\sigma}_1 +\tilde{\sigma}_2
+\tilde{\sigma}_3+\tilde{\sigma}_4),
$$
with
$$
\tilde{\sigma}_1 = \sum_{|n|>N} \sum_{k \neq n} \sum_{j,i\neq n}
\frac{|r(k+j)|^2 |r(i+n)|^2}{|n - k|^2|n - j||n - i|} \leq C \|r\|^2
(\mathcal{E}_{2N} (r))^2 \leq C\|r\|^2 \rho_N^2
$$
(by (\ref{t14}) in Lemma \ref{lemt2});
$$
\tilde{\sigma}_2 = \sum_{|n|>N} \sum_{k \neq n} \sum_{j\neq n}
\frac{|r(k+j)|^2 |r(2n)|^2 }{|n - k|^2|n - j|}
$$
$$
\leq  \sum_{|n|>N} |r(2n)|^2\sum_{k \neq n} \frac{1}{|n-k|^2} \sum_j
|r(k+j)|^2 \leq C \|r\|^2 (\mathcal{E}_{2N} (r))^2 \leq C\|r\|^2
\rho_N^2;
$$
$$
\tilde{\sigma}_3 =\sum_{|n|>N} \sum_{k \neq n}  \sum_{i\neq n}
\frac{|r(k+n)|^2 |r(n+i)|^2}{|n - k|^2|n - i|} $$
$$\leq \sum_{|n|>N} \sum_{k \neq n}
\frac{|r(k+n)|^2}{|n - k|^2} \cdot \sum_i |r(n+i)|^2 \leq C \|r\|^2
\rho_N^2
$$
(by (\ref{t11}) in Lemma \ref{lemt2});
$$
\tilde{\sigma}_4 = \sum_{|n|>N,  k \neq n} \frac{|r(k+n)|^2|r(2n)|^2
}{|n - k|^2} \leq  \|r\|^2 \sum_{|n|>N,  k \neq n} \frac{
|r(k+n)|^2}{|n - k|^2} \leq C \|r\|^2 \rho_N^2
$$
(by (\ref{t11}) in Lemma \ref{lemt2}). These estimates imply the
inequality $\tilde{\sigma} \leq C \|r\|^2 \rho_N^2, $ which
completes the proof of (\ref{148}) for $\nu=2, s > 1.$ \vspace{3mm}

{\em Estimates for $B_3 (s).$} The sums $B_3 (s)$ can be estimated
in a similar way because the indices $k$ and $m$ play symmetric
roles. More precisely, since  $$ B(\lambda,k,i_1, \ldots,i_s,n)
=B(\lambda,n,j_1, \ldots,j_{\tau -1},k)$$ if $j_1 =i_{s -1}, \ldots,
j_{s -1}= i_1,    $  we have $B_3 (s) = B_2 (s).$ Thus, (\ref{148})
holds for $\nu=3.$ \vspace{3mm}

{\em Estimates for $B_4 (s).$} Here $s\geq 1 $ by the definition of
$B_4 (s).$

Fix $s \geq 1$ and consider the sum in (\ref{140}) that defines $B_4
(s);$ then at least one of the indices $j_1, \ldots, j_s$ is equal
to $n.$ Let $\tau \leq t$ be the least integer such that $j_\tau =
n.$ Then, by (\ref{133}) or (\ref{134}), and since $|\lambda
-n|=1/2$ for $\lambda \in C_n,$ we have
$$
B(\lambda,k,j_1,\ldots,j_{\tau -1}, n,j_{\tau +1}, \ldots,j_s,m)=
$$
$$
\frac{1}{2} B(\lambda,k,j_1, \ldots,j_{\tau -1}, n)\cdot B(\lambda,
n,j_{\tau +1}, \ldots,j_s,m).
$$
Therefore,
$$
B_4 (s)  \leq \sum_{\tau=1}^s  \sum_{|n|>N}  \sum_{k \neq
n}\sup_{\lambda \in C_n} \left |\sum_{j_1,\ldots,j_{\tau -1}}
B(\lambda,k,j_1, \ldots,j_{\tau -1},n) \right |^2
$$
$$ \times
 \sum_{m \neq  n} \sup_{\lambda \in C_n} \left |
 \sum_{j_{\tau +1},
\ldots,j_s} B(\lambda,n,j_{\tau +1}, \ldots,j_s,m) \right |^2
$$

On the other hand, by the estimate of $B_3 (s) $ given by
(\ref{148}),
$$
 \sum_{m \neq  n} \sup_{\lambda \in C_n} \left |
 \sum_{j_{\tau +1},
\ldots,j_s} B(\lambda,n,j_{\tau +1}, \ldots,j_s,m) \right |^2 \leq C
\|r\|^2 \rho_N^{2(s-\tau)},\quad |n|>N.
$$
Thus, we have
$$
B_4 (s) \leq  C \|r\|^2 \sum_{\tau=1}^s \rho_N^{2(s-\tau)}
 \sum_{|n|>N}  \sum_{k \neq  n}\sup_{\lambda
\in C_n} \left |\sum_{j_1,\ldots,j_{\tau -1}} B(\lambda,k,j_1,
\ldots,j_{\tau -1},n) \right |^2
$$
Now, by (\ref{148}) for $\nu =2,$
$$
\sum_{|n|>N}  \sum_{k \neq  n}\sup_{\lambda \in C_n} \left
|\sum_{j_1,\ldots,j_{\tau -1}} B(\lambda,k,j_1, \ldots,j_{\tau
-1},n) \right |^2 \leq  C \|r\|^2 \rho_N^{2(\tau -1)}.
$$
Hence,
$$
B_4 (s) \leq  C \|r\|^4 \sum_{\tau=1}^s \rho_N^{2(s-1)}=C s \|r\|^4
 \rho_N^{2(s-1)},
$$
which completes the proof of (\ref{148}).
\end{proof}

Now, we can complete the proof of Theorem \ref{thm1}. Lemma
\ref{lemp}, (\ref{142}) together with the inequalities (\ref{148})
and (\ref{147}) in Proposition \ref{prop1} imply that
\begin{equation}
\label{201} A^{\alpha \beta} (s) \leq 4C \| r\|^2 ( 1 + \|
r\|^2/\rho_N^2 ) (1+s) \rho_N^{2s}
\end{equation}
\begin{equation}
\label{202} \left (A^{\alpha \beta} (s)A^{\alpha \beta} (t) \right
)^{1/2} \leq 4C \| r\|^2 ( 1 + \| r\|^2/\rho_N^2 ) (1+s)(1+t)
\rho_N^{s+t}.
\end{equation}
With $\rho \leq 1/2 $ by (\ref{147}) the inequality (\ref{202})
guarantees that the series on the right side of (\ref{0114})
converges and
$$
\sum_{n>N} \|P_n -P_n^0\|^2 \leq \sum_{n>N} \|P_n -P_n^0\|_{HS}^2
\leq C_1 \| r\|^2 ( 1 + \| r\|^2/\rho_N^2 ) <\infty.
$$
So, Theorem \ref{thm1} is proven subject to Lemmas \ref{lemt1} and
\ref{lemt2} in the next section.
\end{proof}

\section{Technical Lemmas}
In this section we use that
\begin{equation}
\label{t0} \sum_{n>N}\frac{1}{n^2} < \sum_{n>N}\left (
\frac{1}{n-1}- \frac{1}{n}  \right ) = \frac{1}{N}, \quad N \geq 1.
\end{equation}
and
\begin{equation}
\label{t00} \sum_{p \neq \pm n} \frac{1}{(n^2-p^2)^2} <
\frac{4}{n^2}, \quad   n\geq 1.
\end{equation}
Indeed,
$$
\frac{1}{(n^2-p^2)^2} = \frac{1}{4n^2} \left (  \frac{1}{n-p}+
\frac{1}{n+p} \right )^2 \leq \frac{1}{2n^2} \left (
\frac{1}{(n-p)^2}+ \frac{1}{(n+p)^2} \right ).
$$
Therefore, the sum in (\ref{t00}) does not exceed
$$\frac{1}{2n^2}
\left ( \sum_{p \neq \pm n} \frac{1}{(n-p)^2}+\sum_{p \neq \pm n}
\frac{1}{(n+p)^2} \right ) \leq \frac{1}{2n^2}\cdot 2
\frac{\pi^2}{3} <\frac{4}{n^2}.
$$
\vspace{2mm}

\begin{Lemma}
\label{lemt1} If $r = (r(k)) \in \ell^2 (2\mathbb{Z}) $ (or $r =
(r(k)) \in \ell^2 (\mathbb{Z}) $), then
\begin{equation}
\label{t1} \sum_{k\neq n} \frac{|r(n+k)|^2}{|n-k|} \leq
\frac{\|r\|^2}{|n|} + (\mathcal{E}_{|n|} (r))^2, \quad |n| \geq 1;
\end{equation}
\begin{equation}
\label{t2} \sum_{i,k\neq n} \frac{|r(i+k)|^2}{|n-i||n-k|} \leq C
\left ( \frac{\|r\|^2}{\sqrt{n}} + (\mathcal{E}_{|n|} (r))^2 \right
),  \quad |n| \geq 1,
\end{equation}
where $n \in \mathbb{Z},\; i,k \in n+ 2\mathbb{Z} $ (or,
respectively, $i,k \in \mathbb{Z} $) and  $C $ is an absolute
constant.
\end{Lemma}

\begin{proof} If $|n-k| \leq |n|,$ then
we have $|n+k|\geq 2|n|-|n-k|\geq |n|.$ Therefore,
$$
 \sum_{k\neq n}\frac{|r(n+k)|^2}{|n-k|} \leq
\sum_{0<|n-k|\leq |n|} |r(n+k)|^2 +
\sum_{|n-k|>|n|}\frac{|r(n+k)|^2}{|n|}\leq (\mathcal{E}_{|n|} (r))^2
+ \frac{\|r\|^2}{|n|},
$$
which proves (\ref{t1}).

Next we prove (\ref{t2}). We have
\begin{equation}
\label{t5} \sum_{i,k\neq n} \frac{|r(i+k)|^2}{|n-i||n-k|} \leq
\sum_{(i,k) \in J_1} + \sum_{(i,k) \in J_2} + \sum_{(i,k) \in J_3},
\end{equation}
where $ J_1 = \left \{(i,k): \; 0<|n-i|<|n|/2, \; |n-k| < |n|/2
\right \}$,
$$ J_2 = \left \{(i,k): \; i\neq n, \; |n-k|\geq \frac{|n|}{2}
\right \}, \quad J_3 =\left  \{(i,k): \; |n-i|\geq\frac{|n|}{2}, \;
k\neq n  \right \}.
$$

 For $(i,k) \in J_1$ we have $|i+k| = |2n - (n-i)-(n-k)| \geq
2|n|-|n-i|-|n-k|\geq |n|.$ Therefore, by the Cauchy inequality,
$$
\sum_{(i,k) \in J_1} \leq  \left ( \sum_{(i,k) \in J_1}
\frac{|r(i+k)|^2}{|n-i|^2}\right )^{1/2} \left ( \sum_{(i,k) \in
J_1} \frac{|r(i+k)|^2}{|n-k|^2} \right )^{1/2} \leq C
(\mathcal{E}_{|n|} (r))^2.
$$

On the other hand, again by the Cauchy inequality,
$$
\sum_{(i,k) \in J_2}=\sum_{(i,k) \in J_3} \leq \left ( \sum_{(i,k)
\in J_3} \frac{|r(i+k)|^2}{|n-i|^2}\right )^{1/2} \left (
\sum_{(i,k) \in J_3} \frac{|r(i+k)|^2}{|n-k|^2} \right )^{1/2}
$$
$$
\leq \left ( \sum_{|n-i|\geq \frac{|n|}{2}} \frac{1}{|n-i|^2} \sum_k
|r(i+k)|^2 \right )^{1/2} \left( \sum_{k\neq n} \frac{1}{|n-k|^2}
\sum_i |r(i+k)|^2 \right )^{1/2} \leq C\frac{\|r\|^2}{\sqrt{n}},
$$
which completes the proof.
\end{proof}

\begin{Lemma}
\label{lemt2} If $r = (r(k)) \in \ell^2 (2\mathbb{Z}) $ (or $r =
(r(k)) \in \ell^2 (\mathbb{Z}) $), then
\begin{equation}
\label{t11} \sum_{|n|>N,k\neq n} \frac{|r(n+k)|^2}{|n-k|^2} \leq C
\left ( \frac{\|r\|^2}{N} + (\mathcal{E}_N (r))^2 \right );
\end{equation}
\begin{equation}
\label{t12} \sum_{|n|>N}\sum_{i,p \neq n}
\frac{|r(n+i)|^2|r(n+p)|^2}{|n - i||n-p|} \leq  C \left (
\frac{\|r\|^2}{N} + (\mathcal{E}_N (r))^2 \right )\|r\|^2;
\end{equation}
\begin{equation}
\label{t13} \sum_{|n|>N,j, p\neq n}
\frac{|r(j+p)|^2}{|n-j|^2|n-p|^2} \leq C \left ( \frac{\|r\|^2}{N} +
(\mathcal{E}_N (r))^2 \right );
\end{equation}
\begin{equation}
\label{t14} \sum_{|n|>N}\sum_{i,j,p \neq n}
\frac{|r(n+i)|^2|r(j+p)|^2}{|n - i||n-j||n-p|^2} \leq  C \left (
\frac{\|r\|^2}{N} + (\mathcal{E}_N (r))^2 \right )\|r\|^2,
\end{equation}
where $C$ is an absolute constant.
\end{Lemma}

\begin{proof}
With $\tilde{k} = n-k$ and $\tilde{n} = n+k$ it follows that
whenever $|\tilde{k}|\leq |n| $ we have
$|\tilde{n}|=|2n-\tilde{k}|\geq 2|n|-|\tilde{k}|\geq |n|.$
Therefore,
$$
\sum_{|n|>N} \sum_{k\neq n}\frac{|r(n+k)|^2}{|n-k|^2} =
\sum_{|n|>N}\sum_{0<|n-k|\leq |n|} +\sum_{|n|>N}\sum_{|n-k|>|n|}
$$
$$
\leq \sum_{|\tilde{k}|>0} \frac{1}{|\tilde{k}|^2}
\sum_{|\tilde{n}|>N} |r(\tilde{n})|^2 + \sum_{|n|>N}\frac{1}{n^2}
\sum_{k} |r(n+k)|^2 \leq C\left ( (\mathcal{E}_N (r))^2
+\frac{\|r\|^2}{N} \right ),
$$
which proves (\ref{t11}).

Since $\frac{1}{|n-i||n-p| }  \leq \frac{1}{2} \left (
\frac{1}{|n-i|^2} +\frac{1}{|n-p|^2 } \right ),$ the sum in
(\ref{t12}) does not exceed
$$
\frac{1}{2}\sum_{|n|>N,i\neq n} \frac{|r(n+i)|^2}{|n-i|^2} \sum_p
|r(n+p)|^2 + \frac{1}{2}\sum_{|n|>N,p\neq n}
\frac{|r(n+p)|^2}{|n-p|^2} \sum_i |r(n+i)|^2.
$$
In view of (\ref{t11}), the latter is less than $  C \left (
\frac{\|r\|^2}{N} + (\mathcal{E}_N (r))^2 \right )\|r\|^2,$ which
proves (\ref{t12}).

In order to prove (\ref{t13}), we set $\tilde{j}= n-j$ and
$\tilde{p}= n-p. $ Then
$$
\sum_{|n|>N; j,p\neq n} \frac{|r(j+p)|^2}{|n-j|^2|n-p|^2}=
\sum_{\tilde{j},\tilde{p} \neq 0} \frac{1}{\tilde{j}^2}
\frac{1}{\tilde{p}^2}\sum_{|n|>N} |r(2n-\tilde{j}-\tilde{p}|^2
$$
$$
\leq \sum_{0<|\tilde{j}|,|\tilde{p}|\leq N/2} \frac{1}{\tilde{j}^2}
\frac{1}{\tilde{p}^2} \sum_{n>N} |r(2n-\tilde{j}-\tilde{p}|^2
+\sum_{|\tilde{j}| > N/2} \sum_{|\tilde{p}|\neq 0} \cdots +
\sum_{|\tilde{j}|\neq 0} \sum_{|\tilde{p}|> N/2} \cdots
$$
$$
\leq C(\mathcal{E}_N (r))^2 + \frac{C}{N}\|r\|^2
+\frac{C}{N}\|r\|^2,
$$
which completes the proof of (\ref{t13}).

Let $\sigma $ denote the sum in (\ref{t14}). The inequality  $ab
\leq (a^2 + b^2)/2,$ considered with $a=1/|n-i|$ and $b=1/|n-j|,$
implies that $\sigma\leq (\sigma_1 + \sigma_2)/2,$ where
$$
\sigma_1 = \sum_{|n|>N, i\neq n} \frac{|r(n+i)|^2}{|n-i|^2}
\sum_{p\neq n} \frac{1}{|n-p|^2} \sum_j |r(j+p)|^2 \leq C\left (
(\mathcal{E}_N (r))^2 +\frac{\|r\|^2}{N} \right ) \|r\|^2
$$
(by (\ref{t11})), and
$$
\sigma_2 = \sum_{|n|>N} \sum_{j,p\neq n}
\frac{|r(j+p)|^2}{|n-j|^2|n-p|^2} \sum_i |r(n+i)|^2 \leq C\left (
(\mathcal{E}_N (r))^2 +\frac{\|r\|^2}{N} \right )\|r\|^2
$$
(by (\ref{t13})). Thus (\ref{t14}) holds.

\end{proof}

\section{Conclusions}

1. The convergence of the series (\ref{110}) is the analytic core of
Bari--Markus Theorem (see \cite{GK}, Ch.6, Sect.5.3, Theorem 5.2)
which guarantees that the series $\sum_{|n| >N} P_n f $ converges
unconditionally in $L^2$ for every $f \in L^2. $ But in order to
have the identity  $$ f = S_N f +
 \sum_{|n| >N} P_n f ,$$
we need to check  the "algebraic" hypotheses in Bari--Markus
Theorem:

(a)  The system of projections
\begin{eqnarray}
\label{c01} \{S_N; \; \; P_n, \; \; |n| >N \}
\end{eqnarray}
is {\it complete}, i.e., the linear span of the system of subspaces
\begin{eqnarray}
\label{c02} \{E^*; \; \; E_n, \; \; |n| >N \}, \quad E^* = Ran \,
S_N, \; E_n = Ran \, P_n,
\end{eqnarray}
is dense in $L^2 (I).$

(b) The system of subspaces (\ref{c02}) is {\em minimal,} i.e.,
there is no vector in one of these subspaces that belongs to the
closed linear span of all other subspaces. Condition (b) holds
because the projections in (\ref{c01}) are continuous, commute and
$$
P_n S_N =0, \quad P_n P_m =0 \quad \text{for} \; m \neq n, \quad
|m|,|n| >N.
$$

The system (\ref{c01}) is complete; this fact is well known since
the early 1950's (see details in \cite{K,KL,GK}). More general
statements are proven in \cite{MalOr} and \cite{Mit04}, Theorems 6.1
and 6.4 or Proposition 7.1.

Therefore, all hypotheses of Bari--Markus Theorem hold, and we have
the following theorem.

\begin{Theorem}
\label{thm2} Let $L$ be the Dirac operator (\ref{01}) with an
$L^2$-potential $v,$  subject to the boundary conditions $bc =
Per^\pm $ or $Dir.$  Then there is $N \in \mathbb{N}$ such that the
Riesz projections
$$
S_N = \frac{1}{2\pi i} \int_{ |z|= N-1/2 }  (z-L_{bc})^{-1} dz,
\quad P_n = \frac{1}{2\pi i} \int_{ |z-n|= 1/4 }  (z-L_{bc})^{-1} dz
$$
are well--defined, and
 $$ f = S_N f +  \sum_{|n| >N} P_n f, \quad  \forall f \in L^2; $$
moreover, this series converges unconditionally  in $L^2.$
\end{Theorem}

2. General {\em regular} boundary conditions for the operator $L^0$
(or $L$) (\ref{01})--(\ref{02}) are given by a system of two linear
equations
\begin{eqnarray}
\label{c1}  y_1 (0) +b y_1 (\pi) + a y_2 (0) =0
\\ \nonumber d y_1 (\pi) + c y_2 (0) +  y_2 (\pi)=0
\end{eqnarray}
with the restriction
\begin{equation}
\label{c2} bc-ad \neq 0.
\end{equation}
A regular boundary condition is {\em strictly regular,} if
additionally
\begin{equation}
\label{c3}  (b-c)^2 + 4ad \neq 0,
\end{equation}
i.e., the characteristic equation
\begin{equation}
\label{c4} z^2 + (b+c)z + (bc-ad)=0
\end{equation}
has two {\em distinct} roots.

As we noticed in Introduction our main results (Theorem \ref{thm2})
can be extended to the cases of both strictly regular $(SR)$ and
regular  but not strictly regular $(R\setminus SR) \; bc.$ More
precisely, the following  statements hold. \vspace{3mm}

$(SR) \;$ case.   Let $L_{bc}$ be an operator (\ref{01})--(\ref{02})
with $(bc) \in (\ref{c1})-(\ref{c2}).$ Then its spectrum $SP \,
(L_{bc}) = \{\lambda_k, \; k \in \mathbb{Z} \}$ is discrete, $\sup
|Im \,\lambda_k| < \infty, $ $\; |\lambda_k| \to \infty $ as $k \to
\pm \infty,$ and all but finitely many eigenvalues $\lambda_k $  are
simple, $L_{bc} u_k = \lambda_k u_k, \; |k| > N = N(v).$ Put
$$
S_N = \frac{1}{2\pi i } \int_C (z-L_{bc})^{-1} dz,
$$
where the contour $C$ is chosen so that all $\lambda_k, \; |k|\leq
N, $ lie inside of $C, $ and $\lambda_k, \; |k|> N, $ lie outside of
$C. $ Then the spectral decomposition
$$
f= S_N f + \sum_{|k| >  N} c_k (f) u_k, \quad \forall f \in L^2
$$
is well--defined and {\em converges unconditionally} in $L^2.$
\vspace{3mm}

$(R \setminus SR) \;$ case. Let $bc$ be regular, i.e.,
(\ref{c1})-(\ref{c2}) hold, but not strictly regular, i.e.,
\begin{equation}
\label{c5}  (b-c)^2 + 4ad = 0,
\end{equation}
and $z_* = \exp (i\pi \tau)$ be a double root of (\ref{c4}).

Then its spectrum $SP \, (L_{bc}) = \{\lambda_k, \; k \in
\mathbb{Z}\} $ is discrete; it lies in $\Pi_N \cup \bigcup_{m>
N}D_m, \; N=N(v), $ where
$$
\Pi_N = \{z\in \mathbb{C}:  \; |Im \, (z-\tau)|, |Re \, (z-\tau)| <
N-1/2 \}
$$
 and $D_m =\{z\in \mathbb{C}:  \; | (z-m- \tau)|< 1/4  \}.$
The spectral decomposition
$$
f= S_N f + \sum_{|m| >  N} P_m f, \quad \forall f \in L^2
$$
is well--defined if we set
$$
S_N = \frac{1}{2\pi i } \int_{\partial \Pi_N} (z-L_{bc})^{-1} dz,
\quad P_m  = \frac{1}{2\pi i } \int_{\partial D_m}  (z-L_{bc})^{-1}
dz,\quad |m| > N,
$$
and it {\em converges unconditionally} in $L^2.$

Complete presentation and proofs of these general results will be
given elsewhere.

\end{document}